\documentclass[preprint,12pt]{elsarticle}
\usepackage{stmaryrd}

\usepackage{amssymb}
\usepackage{amsthm}
\usepackage{amsmath}
 \newtheorem{thm}{Theorem}[section]
 \newtheorem{cor}[thm]{Corollary}
 \newtheorem{lem}[thm]{Lemma}
 
 \theoremstyle{definition}
 
 \theoremstyle{remark}
 \newtheorem{rem}[thm]{Remark}

\journal{arxiv}

\begin{document}

\begin{frontmatter}



\title{Rigidity Theorem for integral pinched shrinking Ricci solitons}

\author[HPF]{Hai-Ping Fu}
\ead{mathfu@126.com}

\author[LQX]{Li-Qun Xiao}
\ead{xiaoliqun@ncu.edu.cn}

\fntext[fn1]{Supported  by National Natural Science Foundations of China (11261038, 11361041),  Jiangxi Province
Natural Science Foundation of China (20132BAB201005).}

\address[HPF]{Department of Mathematics,  Nanchang University, Nanchang 330031, P.
R. China}

    \address[LQX]{Department of management science and Engineering,  Nanchang University, Nanchang 330031, P.
R. China}

\begin{abstract}
We prove that an $n$-dimensional, $n\geq4$, compact gradient shrinking Ricci
soliton satisfying a $L^{\frac n2}$-pinching condition is isometric to a quotient of the round $\mathbb{S}^n$, which improves the rigidity theorem given by G. Catino\cite{{Ca}}
\end{abstract}

\begin{keyword}
Einstein  manifold\sep Ricci soliton\sep Weyl curvature tensor\sep Yamabe constant

MSC 53C24\sep 53C20
\end{keyword}

\end{frontmatter}


\section{Introduction}
\label{}
In this paper we investigate compact gradient shrinking Ricci solitons satisfying a $L^{\frac n2}$-
pinching condition. Let us recall the concept of Ricci solitons, which was introduced by
Hamilton \cite{Ha} in mid 80's. Let $(M^n, g)$ be an $n$-dimensional,
complete, connected Riemannian manifold. A Ricci soliton is a
Riemannian metric together with a vector field $(M^n,g,X)$ that
satisfies
\begin{equation}
Ric+\frac{1}{2}L_Xg=\lambda g
\end{equation}
for some constant $\lambda$. It is called shrinking, steady or
expanding Ricci soliton depending on whether $\lambda>0$,
$\lambda=0$ or $\lambda<0$ respectively. If there is a smooth
function $f$ on $M$ such that $X=\nabla f$, then the equation (1)
can be written as
\begin{equation}
Ric+\nabla ^2f=\lambda g.\end{equation} This case is called a gradient
(Ricci) soliton. Both
equations (1) and (2) can be considered as perturbations of the
Einstein equation
$$Ric=\lambda g$$ and reduce to this latter case if $X$ or $\nabla
f$ are Killing vector fields. When $X=0$ or $f$ is constant, we call
the underlying Einstein manifold a trivial Ricci soliton.

Ricci solitons are
an important object in the study of the Ricci flow, since they are
self-similar solutions of the flow. They also serve as model cases
of various Harnack inequalities for the Ricci flow, which become
equalities when the flow consists of Ricci solitons. From the
seminal work of  Hamilton \cite{Ha} and Perelman's result \cite{P}
that any compact Ricci soliton is necessarily a gradient soliton, it
is to see that any compact steady or expanding Ricci soliton must be
Einstein \cite{ENM}. So the classification of complete gradient
shrinking solitons plays important roles in the study of the Ricci
flow \cite{C}. In dimension three, T.
Ivey \cite{Iv} proved that the only compact shrinking Ricci solitons are quotients of $\mathbb{S}^3$ with its
standard metric. Dimension four is the lowest dimension where there are interesting examples
of shrinking Ricci solitons.

According to the
decomposition of the Riemannian curvature tensor, for $n\geq 3$, a locally
conformally flat manifold has constant sectional curvature if and only it is Einstein. As a consequence, it follows from the H. Hopf
classification theorem that the space forms are  isometric to the only complete,
simply connected, locally conformally flat, Einstein manifolds. $3$-dimensional connected Einstein manifolds are constant  curvature spaces. Two isolation theorems of Weyl curvature tensor of positive Einstein manifolds are given in \cite{{HV},{IS}}, when its $L^{\frac{n}{2}}$-norm is small.
  The curvature pinching phenomenon plays an important role in global differential
geometry. We are interested in  $L^p$
pinching problems for compact gradient shrinking Ricci solitons.

Now we  introduce the definition of the Yamabe constant. Given a compact  Riemannian $n$-manifold $M$, we consider the Yamabe functional
$$Q_g\colon C^{\infty}_{+}(M)\rightarrow\mathbb{R}\colon f\mapsto Q_g(f)=\frac{\frac{4(n-1)}{n-2}\int_M|\nabla f|^2\mathrm{d}v_g+\int_M Sf^2\mathrm{d}v_g}{(\int_M f^{\frac{2n}{n-2}}\mathrm{d}v_g)^{\frac{n-2}{n}}}.$$
It follows that $Q_g$ is bounded below by H\"{o}lder inequality. We set
$$\mu([g])=\inf\{Q_g(f)|f\in C^{\infty}_{+}(M)\}.$$
This constant $\mu([g])$ is an invariant of the conformal class of $(M, g)$, called the Yamabe constant. 

The important works of Schoen, Trudinger and Yamabe showed that the infimum in the above is always achieved (see \cite{{A},{LP}}). There exists an $f\in C^{\infty}_{+}(M)$ satisfying $Q_g(f)=\mu([g])$. So, the  conformal metric $\tilde{g}=f^{\frac{4}{n-2}}g$ of $g$ has constant scalar curvature $\tilde{S}=\mu([g]){\mathrm{Vol}(\tilde{g})}^{-\frac 2n}$, called the Yamabe metric. Choosing a Yamabe metric in  the  conformal class $[g]$, denoted by $g$, we have the Sobolev inequality
\begin{equation}
\mu([g])(\int_M f^{\frac{2n}{n-2}}\mathrm{d}v_g)^{\frac{n-2}{n}}\leq\frac{4(n-1)}{n-2}\int_M|\nabla f|^2\mathrm{d}v_g+\int_M Sf^2\mathrm{d}v_g, \forall f\in H^2_1(M).
\end{equation}
The Yamabe constant of a given compact manifold
is determined by the sign of scalar curvature \cite{A}. 

In this note, we obtain the following rigidity theorems.
\begin{thm}
Let $(M^n, g)(n\geq4)$ be a complete Einstein $n$-manifold with positive  scalar curvature $S$. For $p\geq \frac n2$,  if
$$\left(\int_{M}|W|^{p}\right)^{\frac 1p}< \left(\frac{\mu([g])}{S}\right)^{\frac {n}{2p}}\left(\frac{2S}{C(n)n}\right),$$
where the constant $C(n)$ is defined in Lemma 2.1,
then $M$ is isometric to a quotient of the round $\mathbb{S}^n$.
\end{thm}
\begin{cor}
Let $(M^n, g)(n\geq4)$ be a complete  Einstein $n$-manifold with positive  scalar curvature $S$. If
$$\left(\int_{M}|W|^{\frac n2}\right)^{\frac 2n}< \frac{2\mu([g])}{C(n)n},$$
then $M$ is isometric to a quotient of the round $\mathbb{S}^n$.
\end{cor}
\begin{rem}
Theorem 1.1 has been proved in \cite{F}. It is easy to see from Corollary 1.2 that the pinching constant is better than the one due to \cite{{Ca},{HV},{IS}}. Theorem 1.1 improves the isolation theorems given by \cite{{Ca},{HV},{IS},{S}}.
\end{rem}

\begin{thm}
Let $(M^n, g)(n\geq 4)$ be an $n$-dimensional compact shrinking Ricci solitions satisfying
\begin{eqnarray}\left(\int_M|W+\frac{\sqrt{2}}{\sqrt{n}(n-2)}\mathring{Ric}\owedge g|^{\frac n2}\right)^{\frac 2n}+\sqrt{\frac{(n-4\gamma)^2(n-1)}{8\gamma^2(n-2)}}\lambda Vol(g)^{\frac n2}\nonumber\\<(\frac{2}{\gamma}-\frac{1}{\gamma^2})\sqrt{\frac{n-2}{32(n-1)}}\mu([g]),\end{eqnarray}
where $\gamma=\frac{n+\sqrt{n^2+8n(n-2)^2}}{8(n-2)}.$ Then $(M^n, g)$ is isometric to a quotient of the round $\mathbb{S}^n$.
\end{thm}
\begin{rem} Integrating equation $(15)$, by the definition of $\mu([g])$, we have 
$$\lambda Vol(g)^{\frac n2}=\frac 1n Vol(g)^{\frac {n-2}{n}}\int_M S\geq \frac 1n\mu([g]).$$
A routine computation gives rise to
\begin{eqnarray*}[(\frac{2}{\gamma}-\frac{1}{\gamma^2})\sqrt{\frac{n-2}{32(n-1)}}+(1-\frac{1}{\gamma})\sqrt{\frac{n-1}{8(n-2)}}]\mu([g])
\geq\sqrt{\frac{n-2}{32(n-1)}}\mu([g])\end{eqnarray*}
and\begin{eqnarray*}
\left(\int_M|W+\frac{\sqrt{2}}{\sqrt{n}(n-2)}\mathring{Ric}\owedge g|^{\frac n2}\right)^{\frac 2n}+\sqrt{\frac{(n-4)^2(n-1)}{8(n-2)}}\lambda Vol(g)^{\frac n2}\\=\left(\int_M|W+\frac{\sqrt{2}}{\sqrt{n}(n-2)}\mathring{Ric}\owedge g|^{\frac n2}\right)^{\frac 2n}+\sqrt{\frac{(n-4\gamma)^2(n-1)}{8\gamma^2(n-2)}}\lambda Vol(g)^{\frac n2}\\+\sqrt{\frac{n^2(\gamma-1)^2(n-1)}{8\gamma^2(n-2)}}\lambda Vol(g)^{\frac n2}\\
\geq\left(\int_M|W+\frac{\sqrt{2}}{\sqrt{n}(n-2)}\mathring{Ric}\owedge g|^{\frac n2}\right)^{\frac 2n}+\sqrt{\frac{(n-4\gamma)^2(n-1)}{8\gamma^2(n-2)}}\lambda Vol(g)^{\frac n2}\\
+(1-\frac{1}{\gamma})\sqrt{\frac{n-1}{8(n-2)}}]\mu([g]).
\end{eqnarray*}
Hence the integral curvature condition in Theorem 1.4, which holds in every dimension $n\geq  4$, is weaker than  the corresponding condition in Theorem 1.4 of \cite{Ca}, holds in every dimension $4\leq n\leq6$. We observe that this result apply also to
complete (possibly non-compact) gradient shrinking Ricci solitons with positive sectional
curvature, since Munteanu and Wang in \cite{MW} recently showed that these conditions force
the manifold to be compact.
\end{rem}

%{Acknowledgement:} The authors would like to thank the referee for
%some helpful suggestions.

\section{Proof of  Theorem 1.1}
\label{}
In what follows, we adopt, without further comment, the moving frame notation with respect to a chosen local orthonormal frame.

Let $M$ be a complete Einstein $n$-manifold. The
decomposition of the Riemannian curvature tensor  into irreducible components yield
\begin{eqnarray*}
R_{ijkl}&=&W_{ijkl}+\frac{1}{n-2}(R_{ik}\delta_{jl}-R_{il}\delta_{jk}+R_{jl}\delta_{ik}-R_{jk}\delta_{il})\nonumber\\
&&-\frac{S}{(n-1)(n-2)}(\delta_{ik}\delta_{jl}-\delta_{il}\delta_{jk})\nonumber\\
&=&W_{ijkl}+\frac{S}{n(n-1)}(\delta_{ik}\delta_{jl}-\delta_{il}\delta_{jk}),
\end{eqnarray*}
where $R_{ijkl}$, $W_{ijkl}$  and $R_{ij}$ denote the components of $Rm$, the Weyl curvature tensor $W$ and the Ricci tensor $Ric$ , respectively,  and  $S$  is the scalar curvature.

Now, we compute the Laplacian of $|W|^2$.
\begin{lem}
Let $M$  be a complete Einstein $n$-manifold with   scalar curvature $S$.   Then
\begin{eqnarray}
\triangle|W|^2\geq2|\nabla W|^2-2C(n)|W|^3+\frac{4S}{n}|W|^2,
\end{eqnarray}
where\begin{equation*}C(n)=
\begin{cases}\frac{\sqrt{6}}{2}, \quad n=4 \\\frac{8\sqrt{10}}{15}, \quad n=5\\\frac{4(n^2+n-4)}{\sqrt{(n-1)n(n+1)(n+2)}}+ \frac{n^2-n-4}{\sqrt{(n-2)(n-1)n(n+1)}}, \quad  n\geq6.
\end{cases}
\end{equation*}
\end{lem}
\begin{rem}  Lemma 2.1 has been proved in \cite{{Ca},{F}}. This estimate is stronger than the one obtained in \cite{{Ca},{F}}.
\end{rem}
\begin{proof}
By the Ricci identities, we obtain
\begin{eqnarray}
\triangle|W|^2&=&2|\nabla W|^2+2\langle W, \triangle W\rangle=2|\nabla W|^2+2W_{ijkl}W_{ijkl,mm}\nonumber\\
&=&2|\nabla W|^2+2W_{ijkl}(W_{ijkm,lm}+W_{ijml,km})\nonumber\\
&=&2|\nabla W|^2+4W_{ijkl}W_{ijkm,lm}\nonumber\\
&=&2|\nabla W|^2+4W_{ijkl}(W_{ijkm,ml}
+W_{hjkm}R_{hilm}\nonumber\\
&&+W_{ihkm}R_{hjlm}+W_{ijhm}R_{hklm}
+W_{ijkh}R_{hmlm})\nonumber\\
&=&2|\nabla W|^2+4W_{ijkl}(W_{hjkm}R_{hilm}
+W_{ihkm}R_{hjlm}\nonumber\\
&&+W_{ijhm}R_{hklm}
+W_{ijkh}R_{hmlm})\nonumber\\
&=&2|\nabla W|^2+4W_{ijkl}(W_{hjkm}W_{hilm}
+W_{ihkm}W_{hjlm}
+W_{ijhm}W_{hklm}
\nonumber\\&&+W_{ijkh}W_{hmlm})
+\frac{4S}{n(n-1)}W_{ijkl}(W_{ljki}+W_{ilkj}+W_{ijlk})
+\frac{4S}{n}|W|^2\nonumber\\
&=&2|\nabla W|^2+4W_{ijkl}(2W_{hjkm}W_{hilm}
-\frac 12W_{ijhm}W_{klhm})+\frac{4S}{n}|W|^2\qquad\nonumber\\
&=&2|\nabla W|^2-4(2W_{ijlk}W_{jhkm}W_{himl}
+\frac 12W_{ijkl}W_{hmij}W_{klhm})+\frac{4S}{n}|W|^2.\qquad
\end{eqnarray}

Case 1. When $n=4$, it was proved in \cite{H} that $$|2W_{ijlk}W_{jhkm}W_{himl}
+\frac 12W_{ijkl}W_{hmij}W_{klhm}|\leq\frac {\sqrt{6}}{4}|W|^3.$$ 

Case 2. When $n=5$,  Jack and  Parker \cite{JP} have proved that $W_{ijkl}W_{hmij}W_{klhm}=4W_{ijlk}W_{jhkm}W_{himl}$. We consider $W$ as a self adjoint operator on $\wedge^2 V$.  Thus by the algebraic inequality for $m$-trace-free symmetric two-tensors $T$, i.e., $tr(T^3)\leq\frac{m-2}{\sqrt{m(m-1)}}|T|^3$, we obtain
\begin{eqnarray*}|2W_{ijlk}W_{jhkm}W_{himl}
+\frac 12W_{ijkl}W_{hmij}W_{klhm}|=|W_{ijkl}W_{hmij}W_{klhm}|\leq\frac {4\sqrt{10}}{15}|W|^3.\end{eqnarray*}

Case 3. When $n\geq 6$, considering $W$ as a self adjoint operator on $S^2 V$, we have
\begin{eqnarray*}|2W_{ijlk}W_{jhkm}W_{himl}
+\frac 12W_{ijkl}W_{hmij}W_{klhm}|\leq2|W_{ijlk}W_{jhkm}W_{himl}|+\frac 12|W_{ijkl}W_{hmij}W_{klhm}|\\
\leq[\frac{2(n^2+n-4)}{\sqrt{(n-1)n(n+1)(n+2)}}+ \frac{n^2-n-4}{2\sqrt{(n-2)(n-1)n(n+1)}}]|W|^3.\end{eqnarray*}

From (6) and Cases 1,2 and 3, we complete the proof of this Lemma.
\end{proof}

By the Kato inequality $|\nabla W|^2\geq \frac{n+1}{n-1}|\nabla |W||^2$ (see \cite{BKN}) and (5), we obtain
\begin{equation} |W|\triangle|W|\geq\frac{2}{n-1}|\nabla |W||^2-C(n)|W|^3+\frac{2S}{n}|W|^2.\end{equation}
Let $u=|W|$.
Using (6), we compute
\begin{eqnarray}
u^{\alpha}\triangle u^{\alpha}&=&u^{\alpha}\left(\alpha(\alpha-1)u^{\alpha-2}|\nabla u|^2+\alpha u^{\alpha-1}\triangle u\right)\nonumber\\
&=&\frac{\alpha-1}{\alpha}|\nabla u^{\alpha}|^2+\alpha
u^{2\alpha-2}u\triangle u\nonumber\\
&\geq&[1-\frac{n-3}{(n-1)\alpha}]|\nabla u^{\alpha}|^2
-C(n)\alpha u^{2\alpha+1}+\frac{2S\alpha }{n} u^{2\alpha},
\end{eqnarray}
where $\alpha$ is a positive constant.

\begin{proof}[{\bf Proof of Theorem
1.1}] When $S>0$, we see from  Myers' Theorem  that $M$ is compact. So Theorem 1.3 in \cite{I} (see also Proposition 3.1 of \cite{LP}) implies that any Einstein metric must be Yamabe, provided it is not conformally flat. We assume that $W$ does not vanish identically. Hence the Einstein metric $g$ is a Yamabe metric in the conformal class $[g]$.

Taking $\alpha=\frac{2p}{n}\geq 1$. Using the Young's inequality, from (8) we obtain
\begin{eqnarray}
u^{\alpha}\triangle u^{\alpha}&\geq&[1-\frac{n-3}{(n-1)\alpha}]|\nabla u^{\alpha}|^2
-{C(n)\epsilon^{1-\alpha}}u^{3\alpha}\nonumber\\
&&-[{C(n)(\alpha-1)\epsilon}-\frac{2S\alpha }{n}]u^{2\alpha}.
\end{eqnarray}
Setting $w=u^{\alpha}$, we can rewrite (9) as
\begin{eqnarray}
w\triangle w&\geq&[1-\frac{n-3}{(n-1)\alpha}]|\nabla w|^2
-{C(n)\epsilon^{1-\alpha}}w^{3}\nonumber\\
&&-[{C(n)(\alpha-1)\epsilon}-\frac{2S\alpha }{n}]w^{2}.
\end{eqnarray}
By (10), we obtain
\begin{eqnarray}
w^{\beta}\triangle w^{\beta}
&\geq&[1-\frac{n-3}{(n-1)\alpha\beta}]|\nabla w^{\beta}|^2
-{C(n)\beta\epsilon^{1-\alpha}}w^{2\beta+1}\nonumber\\
&&-\beta[{C(n)(\alpha-1)\epsilon}-\frac{2S\alpha }{n}]w^{2\beta},
\end{eqnarray}
where $\beta$ is a positive constant. Integrating by parts over $M$,  we get
\begin{eqnarray}
[2-\frac{n-3}{(n-1)\alpha\beta}]\int_{M}|\nabla w^\beta|^2-C(n)\beta\epsilon^{1-\alpha}\int_{M}w^{2\beta+1}\nonumber\\
-\beta[{C(n)(\alpha-1)\epsilon}-\frac{2S\alpha }{n}]\int_{M}w^{2\beta}\leq 0.
\end{eqnarray}
By the H\"{o}lder inequality and (12), we have
\begin{eqnarray}
[2-\frac{n-3}{(n-1)\alpha\beta}]\int_{M}|\nabla w^\beta|^2-C(n)\beta\epsilon^{1-\alpha}(\int_{M}w^{\frac{2n\beta}{n-2}})^{\frac{n-2}{n}}(\int_{M}w^{\frac n2})^{\frac{2}{n}}\nonumber\\
-\beta[{C(n)(\alpha-1)\epsilon}-\frac{2S\alpha }{n}]\int_{M}w^{2\beta}\leq 0.
\end{eqnarray}
Set $\frac {1}{\alpha\beta}=\frac{(n-1)(1+\sqrt{1-\frac{8(n-3)}{n(n-2)\alpha}})}{n-3}$ and $\epsilon=\frac{2S}{C(n)n}$. Combining (3) with (13), we get
\begin{eqnarray}
\left[\left(2-\frac{(n-3)}{(n-1)\alpha\beta}\right)\frac{(n-2)S}{4(n-1)}
-C(n)\beta\epsilon^{1-\alpha}\left(\int_{M}|W|^{p}\right)^{\frac{2}{n}}\right]\left(\int_{M}w^{\frac{2n\beta}{n-2}}\right)^{\frac{n-2}{n}}\leq 0.
\end{eqnarray}
We choose $\left(\int_{M}|W|^{p}\right)^{\frac 1p}< \left(\frac{\mu([g])}{S}\right)^{\frac {n}{2p}}\left(\frac{2S}{C(n)n}\right)$ such that (14) implies $\left(\int_{M}w^{\frac{2n\beta}{n-2}}\right)^{\frac{n-2}{n}}=0$, that is, $W=0$, i.e., $M$ is  Einstein manifold and locally conformally flat manifold. Hence $M$ is isometric to a quotient of the round $\mathbb{S}^n$.
\end{proof}

\section{Proof of  Theorem 1.4}
\label{}
 First, we recall the
following well known formulas (for the proof see [9]).
\begin{lem}
Let $(M^n, g)$ be a  gradient
 Ricci solitons, then the following formulas
hold,
\begin{eqnarray}
\triangle f=n\lambda-S,\\
\triangle_f S=2\lambda S-2|Ric|^2,\\
\triangle_f R_{ik}=\frac {2}{(n-2)(n-1)}\left(S^2\delta_{ik}-nSR_{ik}+2(n-1)R_{ij}R^{j}_{k}-(n-1)|Ric|^2\delta_{ik}\right)\nonumber\\
+2\lambda R_{ik}-2W_{ijkl}R_{jl},
\end{eqnarray}
where the $\triangle_f$ denotes the $f$-Laplacian, $\triangle_f=\triangle-\nabla_{\nabla f}$.
\end{lem}
Integrating equation (15), one has $n\lambda=\frac{\int_M S}{Vol(g)}\geq S_{min}$. When $S$ gets its minimum, from $(16)$, we get
$$\triangle S_{min}\leq\frac{2S_{min}}{n}(n\lambda-S_{min}).$$
This relation, by the strong maximum principle, implies that if $S$ is nonconstant, then
it must be positive everywhere, hence  $\lambda$ and $\mu([g])$ are positive.

\begin{proof}[{\bf Proof of Theorem
1.4}]
A simple computation shows the following equation for the $f$-Laplacian of the squared norm
of the treceless Ricci tensor.
\begin{eqnarray}\frac 12\triangle_f|\mathring{Ric}|^2=\frac 12\triangle|\mathring{Ric}|^2-\frac 12\langle\nabla f, \nabla|\mathring{Ric}|^2\rangle=|\nabla\mathring{Ric}|^2+2\lambda|\mathring{Ric}|^2-2W_{ijkl}\mathring{R}_{ik}\mathring{R}_{jl}\nonumber\\+\frac {4}{n-2}\mathring{R}_{ij}\mathring{R}_{jk}\mathring{R}_{ki}-\frac {2(n-2)}{n(n-1)}S|\mathring{Ric}|^2.\end{eqnarray}
Using Kato inequality, i.e., $|\nabla|\mathring{Ric}||^2\leq |\nabla\mathring{Ric}|^2$ at every point where $|\mathring{Ric}|\neq 0$, we have
\begin{eqnarray}|\mathring{Ric}|\triangle|\mathring{Ric}|\geq2\lambda|\mathring{Ric}|^2-2W_{ijkl}\mathring{R}_{ik}\mathring{R}_{jl}+\frac {4}{n-2}\mathring{R}_{ij}\mathring{R}_{jk}\mathring{R}_{ki}\nonumber\\-\frac {2(n-2)}{n(n-1)}S|\mathring{Ric}|^2+\frac 12\langle\nabla f, \nabla|\mathring{Ric}|^2\rangle.\end{eqnarray}

Set $u=|\mathring{Ric}|$.
By (19) and Proposition 2.1 in \cite{Ca} , we compute
\begin{eqnarray}
u^{\gamma}\triangle u^{\gamma}&=&u^{\gamma}\left(\gamma(\gamma-1)u^{\gamma-2}|\nabla u|^2+\gamma u^{\gamma-1}\triangle u\right)\nonumber\\
&=&\frac{\gamma-1}{\gamma}|\nabla u^{\gamma}|^2+\gamma
u^{2\gamma-2}u\triangle u\nonumber\\
&\geq&(1-\frac{1}{\gamma})|\nabla u^{\gamma}|^2+2\gamma\lambda u^{2\gamma}+\gamma(-W_{ijkl}\mathring{R}_{ik}\mathring{R}_{jl}+\frac {2}{n-2}\mathring{R}_{ij}\mathring{R}_{jk}\mathring{R}_{ki})u^{2\gamma}\nonumber\\
&&-\frac {2(n-2)}{n(n-1)}\gamma Su^{2\gamma}
+\frac 12\gamma u^{2\gamma-2}\langle\nabla f, |\mathring{Ric}|^2\nabla u^2\rangle\nonumber\\
&\geq&(1-\frac{1}{\gamma})|\nabla u^{\gamma}|^2+2\gamma\lambda u^{2\gamma}-\sqrt{\frac {2(n-2)}{n-1}}\gamma(|W|^2+\frac{8}{n(n-2)}u^2)^{\frac 12}u^{2\gamma}\nonumber\\
&&-\frac {2(n-2)}{n(n-1)}\gamma Su^{2\gamma}+\frac 12\langle\nabla f, \nabla u^{2\gamma}\rangle.
\end{eqnarray}
Integrating by parts over $M^n$ and using equation (15) it follows that
\begin{eqnarray}
0\geq-\frac 12 \int_M u^{2\gamma}\triangle f+(2-\frac{1}{\gamma})\int_M |\nabla u^{\gamma}|^2+2\gamma\lambda\int_M u^{2\gamma}\nonumber\\
-\sqrt{\frac {2(n-2)}{n-1}}\gamma\int_M (|W|^2+\frac{8}{n(n-2)}u^2)^{\frac 12}u^{2\gamma}-\frac {2(n-2)}{n(n-1)}\gamma \int_M Su^{2\gamma}\nonumber\\
\geq(2-\frac{1}{\gamma})\int_M |\nabla u^{\gamma}|^2-\sqrt{\frac {2(n-2)}{n-1}}\gamma\int_M (|W|^2+\frac{8}{n(n-2)}u^2)^{\frac 12}u^{2\gamma}\nonumber\\-\frac{n-4\gamma}{2}\lambda\int_M u^{2\gamma}
+\frac {n(n-1)-4(n-2)\gamma }{2n(n-1)}\int_M Su^{2\gamma}.
\end{eqnarray}
For $2-\frac{1}{\gamma}>0$, by the definition of Yamabe constant and (21), we get
\begin{eqnarray}
0\geq(2-\frac{1}{\gamma})\frac{n-2}{4(n-1)}\mu([g])\left(\int_M  u^{\frac{2n\gamma}{n-2}}\right)^{\frac{n-2}{n}}-\sqrt{\frac {2(n-2)}{n-1}}\gamma\int_M (|W|^2+\frac{8}{n(n-2)}u^2)^{\frac 12}u^{2\gamma}\nonumber\\-\frac{n-4\gamma}{2}\lambda\int_M u^{2\gamma}
+\frac {2n+\frac{1}{\gamma}n(n-2)-8(n-2)\gamma }{2n(n-1)}\int_M Su^{2\gamma}.
\end{eqnarray}
By H\"{o}lder inequality, since $\lambda>0$ and $n-4\gamma\geq0$, we obtain
\begin{eqnarray*}
0\geq\left[(2-\frac{1}{\gamma})\frac{n-2}{4(n-1)}\mu([g])-\sqrt{\frac {2(n-2)}{n-1}}\gamma\left(\int_M (|W|^2+\frac{8}{n(n-2)}u^2)^{\frac n4}\right)^{\frac 2n}-\frac{n-4\gamma}{2}\lambda Vol(g)^{\frac 2n}\right]\nonumber\\\left(\int_M  u^{\frac{2n\gamma}{n-2}}\right)^{\frac{n-2}{n}}
+\frac{2n+\frac{1}{\gamma}n(n-2)-8(n-2)\gamma}{2n(n-1)}\int_M Su^{2\gamma}.
\end{eqnarray*}
Since $W$ is totally trace-free, one has
$$|W+\frac{\sqrt{2}}{\sqrt{n}(n-2)}\mathring{Ric}\owedge g|^2=|W|^2+\frac{8}{n(n-2)}|\mathring{Ric}|^2$$
and the pinching condition (4) implies that $(M^n, g)$ is Einstein. We assume that $W$ does not vanish identically. Hence the Einstein metric $g$ is a Yamabe metric in the conformal class $[g]$, i.e.,
$$\mu([g])=Vol(g)^{\frac {n-2}{n}}\int_M S.$$
Moreover, integrating equation (15) one has
$$\lambda Vol(g)^{\frac n2}=\frac 1n Vol(g)^{\frac {n-2}{n}}\int_M S=\frac 1n\mu([g]).$$
Hence, the pinching condition (4) implies
\begin{equation}\left(\int_M |W|^{\frac n2}\right)^{\frac 2n}<\frac {2}{n\sqrt{2(n-2)(n-1)}}\mu([g]).\end{equation}
It is sufficient to observe
that the pinching constant in (23) is strictly smaller than the one in Corollary 1.2. This completes the proof of Theorem 1.4.
\end{proof}


\begin{thebibliography}{00}
\bibitem{A}
T. Aubin,  \textit{Some Nonlinear Problems in Riemannian Geometry.} Springer-Verlag, Berlin
1998.
\bibitem{BKN}
S. Bando, A. Kasue and H. Nakajima, \textit{On a construction at infinity on manifolds with fast curvature decay and maximal volume growth.} Invent. Math. \textbf{97} (1989), 313--349.

\bibitem{B}
A. L. Besse, \textit{Einstein manifolds.} Springer-Verlag, Berlin, 1987.


\bibitem{C}
H. D. Cao, \textit{Recent progress on Ricci solitons}.arXiv:0908.2006v1 [math.DG], 2009.


\bibitem{Ca}
G. Cation, \textit{Integral pinched shrinking Ricci solitons.} arXiv:1509.07416vl [math.DG], 2015.


\bibitem{ENM}
 M. Eminenti, G. La Nave and C. Mantegazza, \textit{Ricci solitons: the equation point of view.}
 Manuscript Math. \textbf{127} 2008,  345--367.
\bibitem{F}
H. P. Fu and L. Q. Xiao, \textit{Einstein  manifolds with finite  $L^p$-norm of the Weyl curvature.} Submitted to Differ.
Geom. Appl. in September 1st 2015.


\bibitem{Ha}
R. S. Hamilton, \textit{The Ricci flow on surfaces}. Mathematics and general relativity, (Santa Cruz,CA, 1986), volume 71 of Contemp. Math., pages 237--262. Am. Math. Soc., 1988.

\bibitem{HV}
E. Hebey  and  M. Vaugon,  \textit{Effective $L^p$ pinching for the concircular curvature.} J. Geom. Anal. \textbf{6}
(1996),  531--553.
\bibitem{H}
G. Huisken,  \textit{Ricci deformation of the metric on a Riemannian manifold.} J. Differential Geom. \textbf{21} (1985),
47--62.
\bibitem{I}
M. Itoh,  \textit{Yamabe metrics and the space of conformal structures.} Intern.
J. Math. \textbf{2} (1991), 659--671.

\bibitem{IS}
M. Itoh and H. Satoh,  \textit{Isolation of the Weyl conformal tensor for Einstein manifolds.} Proc.
Jpn. Acad. A \textbf{78} (2002), 140--142.
\bibitem{Iv}
T. Ivey,  \textit{Ricci solitons on compact three¨Cmanifolds.} Differential Geom. Appl. \textbf{3} (1993), 301--307
\bibitem{JP}
I. Jack and L. Parker, \textit{Linear independence of renormalisation counterterms in curved space-times of
arbitrary dimensionality.} J. Math. Phys. \textbf{28} (1987), 1137--1139.


\bibitem{LP}
M. J. Lee and T. H. Parker,  \textit{The Yamabe problem.} Bull.
A. M. S. \textbf{17} (1987), 37--91.
\bibitem{MW}
O. Munteanu and J. Wang, \textit{Positively curved shrinking Ricci solitons are compact.} arXiv: 1504.07898vl
[math.DG],  2015.

\bibitem{P}
G. Perelman, \textit{The entropy formula for the Ricci flow and its
geometric applications}. arXiv:math/0211159vl[math.DG], 2002.



\bibitem{S}
M. Singer,  \textit{Positive Einstein metrics with small $L_{n/2}$-norm of the Weyl tensor.} Differ.
Geom. Appl. \textbf{2} (1992), 269--274.





\end{thebibliography}
\end{document}